\documentclass[11pt,a4paper]{article}
\usepackage{amssymb}

\usepackage{graphicx}
\usepackage{amsmath}
\usepackage{ams}

\def\mineappendix{
        \setcounter{section}{1}
        \setcounter{subsection}{0}
        \def\thesection{\Alph{section}}
        \def\sectionap{\@startsection  {section}{1}{\z@}
                        {-3.5ex plus-1ex minus-.2ex} {0ex plus.2ex}
                        {\reset@font\Large\bf  Appendix:  \, }}}
\makeatother
\def\Proclaim #1. #2\par{\bigbreak\noindent{\sc#1.\enspace}{\it#2}\par}

\font\Bbbfont=msbm10
\newfam\msbfam
\textfont\msbfam=\Bbbfont \textfont\msbfam=\Bbbfont

\def\pa{\partial}

\def\la{\lambda}
\def\Om{\Omega}

\def\om{\omega}
\def\Om{\Omega}

\newtheorem{Theorem}{Theorem}
\newtheorem{Lemma}{Lemma}
\newtheorem{Corollary}{Corollary}

\usepackage[dvipdfm, pdfstartview=FitH,
CJKbookmarks=true, bookmarksnumbered=true,
bookmarksopen=true,
citecolor=blue 
,colorlinks=true
,pdfborder=1 
,pdfstartpage=2 
,pdfstartview=Fit]{hyperref} 
\allowdisplaybreaks
\title{\bf An Isoperimetric Inequality for Eigenvalues of the Bi-Harmonic Operator}
\author{Qing Ding,\footnote{Email:qding@fudan.edu.cn; Fax: 86-21-65646073.}
\ \ \ Feng Gao\ \ \ and\ \ \   Yongqian Zhang\footnote{Email:yongqianz@fudan.edu.cn.}   \\
\it School of Mathematical Sciences and Key Lab. of Math. for Nonlinear Sciences  \\
\it Fudan University, Shanghai 200433, P.R. China\\}
\date{}
\begin{document}
\maketitle

\begin{abstract}
In this article, we put forward a Neumann eigenvalue problem for the
bi-harmonic operator $\Delta^2$ on a bounded smooth domain $\Om$ in
the Euclidean $n$-space ${\bf R}^n$ ($n\ge2$) and then prove that
the corresponding first non-zero eigenvalue $\Upsilon_1(\Om)$ admits
the isoperimetric inequality of Szeg\"o-Weinberger type:
$\Upsilon_1(\Om)\le \Upsilon_1(B_{\Om})$, where $B_{\Om}$ is a ball
in ${\bf R}^n$  with the same volume of $\Om$. The isoperimetric
inequality of Szeg\"o-Weinberger type for the first nonzero Neumann
eigenvalue of the even-multi-Laplacian operators $\Delta^{2m}$
($m\ge1$) on $\Om$ is also exploited.
\end{abstract}

Keywords: Neumann problem; first eigenvalue; Min-Max principle;
trial function

\smallskip
Mathematics Subject Classification 2000: 58J50, 35P15, 74K20, 35J40,
35R01

\section{\large Introduction}
The study of eigenvalue problems of the Laplacian operator $\Delta$
and related operators on a bounded smooth domain $\Om$ in the
Euclidean space ${\bf R}^n$ or in an $n$-dimensional Riemannian
manifold $M$ is a topic that always attracts much attentions. The
motivation comes directly from physical problems such as the
membrane problem, the vibration of a clamped plate, the buckling
problem and so on (for example, see \cite{CH}). The typical
eigenvalue problems for the Laplacian operator $\Delta$ are
\begin{eqnarray}
\left\{\begin{array}{c} \Delta u=-\lambda u, ~~~\hbox{in}~~\Om,
\\
u=0,~~~~~~~\hbox{on}~~\pa\Om\end{array}\right. ~~~~~\hbox{(Dirichlet
Problem)}
 \label{D}
\end{eqnarray}
and
\begin{eqnarray}
\left\{\begin{array}{c} \Delta u=-\mu u, ~~\hbox{in}~~\Om,
\\
\frac{\pa u}{\pa n}=0, ~~~~~\hbox{on}~~\pa\Om,\end{array}\right.
~~~~~\hbox{(Neumann Problem)}
 \label{N}
\end{eqnarray}
where $\frac{\pa}{\pa n}$ denotes the outward normal derivative on
the boundary $\pa\Om$ of $\Om$. We always assume $\pa\Om$ is
sufficiently smooth throughout the article, unless otherwise
specified. One of the most important and interesting subjects in the
study of eigenvalue problems is to understand the isoperimetric
property of eigenvalues (see \cite{A,Ch,Ch1}, for example). The
earliest isoperimetric inequality for eigenvalues is the first
Dirichlet eigenvalue of the Laplacian $\Delta$ (i.e. the fixed
membrane) conjectured by Rayleigh in 1877 (see \cite{Ra}):
\begin{eqnarray}
\la_1(\Om)\ge \la_1(B_{\Om}),\label{Ray}
\end{eqnarray}
 for a bounded domain
$\Om\subset {\bf R}^n$, where $B_{\Om}$ is a ball in ${\bf R}^n$
with the same volume of $\Om$,  and the equality holds if and only
if $\Om$ is the ball $B_{\Om}$. This conjecture was proved
independently by Faber \cite{Fa} and Krahn \cite{Kr} by using the
technique of symmetrization. The second isoperimetric inequality for
eigenvalues is
\begin{eqnarray}
\mu_1(\Om)\le \mu_1(B_{\Om})\label{Wei}
\end{eqnarray}
for the first nonzero Neumann eigenvalue of the Laplacian $\Delta$
(i.e. the free membrane) on a bounded smooth domain $\Om$ in ${\bf
R}^n$ and the equality holds if and only if $\Om$ is the ball
$B_{\Om}$. This isoperimetric inequality was proved in the case of
$n=2$ and conjectured for any $n\ge3$ by Szeg\"o in \cite{Se} in
1954. Two years later, Weinberger showed that the Szeg\"o's
conjecture is true in \cite{Wei}. A beautiful isoperimetric
inequality for eigenvalues was conjectured by Payne, Polya and
Weiberger in 1956 in \cite{PPW}, which says
\begin{eqnarray}
\frac{\la_2(\Om)}{\la_1(\Om)}\le\frac{\la_2(B_{\Om})}{\la_1(B_{\Om})}\label{AB}
\end{eqnarray}
for the first two eigenvalues of the Dirichlet Laplacian of the
Laplacian $\Delta$ on a bounded smooth domain $\Om$ in ${\bf R}^n$
and the equality holds if and only if $\Om$ is the ball $B_{\Om}$.
This conjecture was proved by Ashbaugh and Benguria in \cite{AB,AB1}
in 1991. One may also refer to \cite{A} and \cite{Ch} Chapter IV for
the discussion of the isoperimetric property of the first
(Dirichlet) eigenvalue of the Laplacian operator $\Delta$ on a
normal domain in a general Riemannian manifold $M$, especially in a
space form ${\mathbb M}_{\kappa}$ of constant curvature $\kappa$.

On the other hand, the following Dirichlet eigenvalue problem of the
bi-harmonic operator $\Delta^2$ on a bounded domain $\Om$ in ${\bf
R}^n$ is well-known:
\begin{eqnarray}
\left\{\begin{array}{c} \Delta^2 u=\Gamma u, ~~~\hbox{in}~~\Om,
\\
u=\frac{\pa u}{\pa n}=0,~\hbox{on}~\pa\Om.\end{array}\right.
 \label{3}
\end{eqnarray}
When $n=2$, Problem (\ref{3})  relates to the vibration of a clamped
plate in physics (see \cite{CH}). It has discrete eigenvalues
$$0<\Gamma_1\le\Gamma_2\le\cdots\le \Gamma_k\le \cdots\nearrow \infty.$$
In 1877, Rayleigh also made a conjecture that, for the first
eigenvalue $\Gamma_1(\Om)$ of (\ref{3}) on $\Om$ in the plane, one
has the following isoperimetric inequality
\begin{eqnarray}
\Gamma_1(\Om)\ge\Gamma_1(B_{\Om}),\label{Bi}
\end{eqnarray}
and the equality holds if and only if $\Om$ is a disk $B_{\Om}$.
This isoperimetric inequality might apply equally well on a bounded
smooth domain $\Om$ in the Euclidean $n$-space ${\bf R}^n$
($n\ge3$), which is regarded as the inequality of
Rayleigh-Faber-Krahn type for the first Dirichlet eigenvalue of the
bi-harmonic operator $\Delta^2$. Conjecture (\ref{Bi}) was proved by
Nadirashvili in the case of $n=2$ in \cite{Na} and later by Ashbaugh
and Benguria in dimension $n=3$ in \cite{AB2}. However, this
conjecture is still open for $n\ge4$.

It is very natural to conjecture that the isoperimetric inequality
of Szeg\"o-Weinberger type might hold equally well for the
bi-harmonic operator $\Delta^2$ on a bounded smooth domain $\Om$ in
the Euclidean $n$-space ${\bf R}^n$. In this article, we shall focus
our attention on this problem. The understanding of this problem is
in fact divided into two parts. One is to propose a suitable
eigenvalue problem for the bi-harmonic operator $\Delta^2$ on a
bounded smooth domain $\Om$ in ${\bf R}^n$, which is actually a
generalization the Neumann eigenvalue problem (\ref{N}) for the
bi-harmonic $\Delta^2$. Another one is to prove that the
isoperimetric inequality of Szeg\"o-Weinberger type is true for the
corresponding first nonzero eigenvalue.

The first problem is not difficult to answer. In fact, by using the
Green's formula on a bounded smooth domain $\Om$ in ${\bf R}^n$:
$$
\int_{\Om}(u\Delta v-\Delta u v)=\oint_{\pa \Om}\left(u\frac{\pa
v}{\pa n}-v\frac{\pa u}{\pa n}\right),$$ where, as mentioned,
$\frac{\pa}{\pa n}$ denotes the outward normal derivative on the
boundary $\pa\Om$ of $\Om$, we have
\begin{eqnarray}
\int_{\Om}\Delta u\cdot\Delta v=\int_{\Om}u\cdot \Delta^2
v-\oint_{\pa\Om}\left(u\frac{\pa\Delta v}{\pa n}-\frac{\pa u}{\pa
n}\Delta v\right)\label{1}
\end{eqnarray}
and
\begin{eqnarray}
\int_{\Om}\Delta u\cdot\Delta v=\int_{\Om}\Delta^2u\cdot
v-\oint_{\pa\Om}\left(v\frac{\pa\Delta u}{\pa n}-\frac{\pa v}{\pa
n}\Delta u\right).\label{2}
\end{eqnarray}
Therefore, besides the Dirichlet eigenvalue problem (\ref{3}), the
formulas (\ref{1}) and (\ref{2}) lead us to propose the following
three possible eigenvalue problems for the bi-harmonic operator
$\Delta^2$ such that, under the boundary conditions, $\Delta^2$ is a
self-adjoint operator from $L^2(\Om)$ to $L^2(\Om)$:
\begin{eqnarray}
\left\{\begin{array}{c} \Delta^2 u=\Upsilon u, ~~~~~~\hbox{in}~~\Om,
\\
\frac{\pa u}{\pa n}=\frac{\pa\Delta u}{\pa
n}=0,~\hbox{on}~\pa\Om,\end{array}\right.
 \label{4}
\end{eqnarray}
\begin{eqnarray}
\left\{\begin{array}{c} \Delta^2 u=\Gamma u, ~~~~~~~\hbox{in}~~\Om,
\\
u=\Delta u=0,~\hbox{on}~\pa\Om,\end{array}\right.\label{401}
\end{eqnarray}
and
\begin{eqnarray}
\left\{\begin{array}{c} \Delta^2 u=\Upsilon u,
~~~~~~~~~\hbox{in}~~\Om,
\\
\Delta u=\frac{\pa\Delta u}{\pa
n}=0,~\hbox{on}~\pa\Om.\end{array}\right.
 \label{402}
\end{eqnarray}
Among the candidates problems (\ref{4}-\ref{402}), only Problem
(\ref{4}) involves the usual Neumann condition: $\frac{\pa u}{\pa
n}|_{\pa\Om}=0$. Furthermore, with the boundary conditions proposed
by (\ref{4}), it is easy to verify that the bi-harmonic operator
$\Delta^2$ is a self-adjoint elliptic differential operator from
$L^2(\Om)$ to $L^2(\Om)$ with densely defined domain
$${\cal D}(\Delta^2)=\left\{u\in H^4(\Om)~:~\frac{\pa u}{\pa n}\bigg|_{\pa\Om}
=\frac{\pa\Delta u}{\pa n}\bigg|_{\pa\Om}=0\right\}. $$ From the
standard spectral theory in functional analysis, the eigenvalue
problems (\ref{4}) has discrete eigenvalues
$$0=\Upsilon_0<\Upsilon_1\le\Upsilon_2\le\cdots\le \Upsilon_k\le \cdots\nearrow \infty.$$
Therefore, we settle naturally Problem (\ref{4}) as the
generalization the Neumann eigenvalue problem (\ref{N}) for the
bi-harmonic $\Delta^2$. What remains for us to do is now to show
that the first non-zero eigenvalue $\Upsilon_1(\Om)$ of the
bi-harmonic operator $\Delta^2$ given by the problem (\ref{4})
admits the isoperimetric property:
\begin{eqnarray}
\Upsilon_1(\Om)\le \Upsilon_1(B_{\Om}),\label{?}
\end{eqnarray}
where $B_{\Om}$ is a ball in ${\bf R}^n$ with the same volume of
$\Om$, and the equality holds if and only if $\Om$ itself is a ball
$B_{\Om}$.

The aim of this paper is to prove that the isoperimetric inequality
(\ref{?}) does hold for any bounded smooth domain $\Om$ in the
Euclidean space ${\bf R}^n$ ($n\ge 2$) and then to give an
affirmative answer to our problem mentioned above. The method
applied here is mainly the Min-Max principle for eigenvalues and a
suitably choice of trial functions. We also display the
isoperimetric inequality of Szeg\"o-Weinberger type for the first
nonzero Neumann eigenvalue of the even-multi-Laplacian operators
$\Delta^{2m}$ ($m\ge2$) on $\Om$ in ${\bf R}^n$ ($n\ge2$). Up to the
authors' knowledge, there are no other results in literature on
eigenvalues of the multi-Laplacian operators than universal
inequalities (see \cite{CY,WX} and the references there in). Though
we believe that the similar isoperimetric inequality (\ref{?})
should hold for the first nonzero Neumann eigenvalue of the
odd-multi-Laplacian $\Delta^{2m+1}$ ($m\ge1$) on a bounded domain
$\Om$ in ${\bf R}^n$, our method applied in this article does not
work directly in this case.

The paper is organized as follows. In section 2, we show that the
first nonzero Neumann eigenvalue of the bi-harmonic operator
$\Delta^2$ on a ball $B_R(0)$ with radius $R>0$ in ${\bf R}^n$
($n\ge2$) is exactly the square of the first Neumann nonzero
eigenvalue of the Laplacian $\Delta$ on the ball. The $n$ first
Neumann eigenfunctions of $\Delta$ on $B_R(0)$ in ${\bf R}^n$ are
shown to be the first Neumann eigenfunctions of $\Delta^2$ on the
ball. The same conclusion are proved to be hold for the
even-multi-Lapacian $\Delta^{2m}$ with $m\ge1$. In section 3, we
prove that the isoperimetric inequality of Szeg\"o-Weinberger type
for the first nonzero Neumann eigenvalue $\Upsilon_1(\Om)$ of
$\Delta^2$ is true on any bounded smooth domain $\Om$ in ${\bf R}^n$
($n\ge2$). The similar isoperimetric inequality for the first
nonzero Neumann eigenvalue of $\Delta^{2m}$ is also shown to be
true. Some remarks and questions are given in this section.

\section{\large Basic lemmas}
A key step in proving the isoperimetric inequality (\ref{?}) is to
well understand the first nonzero Neumann eigenvalue of the
bi-harmonic operator $\Delta^2$ and its eigenfunctions on a ball in
the Euclidean $n$-space ${\bf R}^n$. One may refer to \cite{CH} for
the characterization of the Dirichlet eigenvalues of the bi-harmonic
operator $\Delta^2$ and their eigenfunctions on a disk in the plane.
In this section, we shall show that the first nonzero Neumann
eigenvalue $\Upsilon_1(B_R(0))$ of the bi-harmonic operator
$\Delta^2$ on the ball $B_R(0)$ with the center at the origin and
radius $R$ in ${\bf R}^n$ is $\mu_1^2(R)$, where $\mu_1(R)$ is the
first nonzero Neumann eigenvalue of the Laplacian $\Delta$ on
$B_R(0)$.

Let $\{g_{k}(\om), k\ge 0\}$ be the set of all the eigenfunctions of
$\Delta_{S^{n-1}}$ on the $(n-1)$-sphere $S^{n-1}$, where $\om$
denotes the spherical coordinates in $S^{n-1}$. We known from
\cite{Ch,CH} that
$$
g_0(\om)=\hbox{constant},~g_k(\om)=x_k\bigg|_{S^{n-1}},~(k=1,\cdots,n),
\cdots
$$
where $(x_1,\cdots,x_n)$ is the standard Euclidean coordinates on
${\bf R}^n$. Furthermore $\{g_{k}(\om), k\ge 0\}$ forms an
orthogonal basis in $L^2(S^{n-1})$ up to a renormalization, that is,
\[\Delta_{S^{n-1}}g_{k}=\lambda_k g_{k},~k\ge0 \]
with \[ \int_{S^{n-1}} |g_k(\omega)|^2 d\omega =1,\] where
$\{\lambda_k\}_{k=0}^{\infty}$ are eigenvalues of the $(n-1)$-sphere
$S^{n-1}$ counted without multiplicity.  From the standard spectral
theory of spheres (see \cite{Ch,CH} also), we know that
$\lambda_k=-j(j+n-2)$ for some $j\ge 0$. The first few of them are
$$
\lambda_0=0,~\lambda_1=\cdots=\lambda_n=-(n-1).$$ For a bounded
smooth domain $\Om$, since $\Delta^2$ is a densely defined
self-adjoint elliptic differential operator from $L^2(\Om)$ to
$L^2(\Om)$, we see that the first nonzero Neumann eigenvalue
(\ref{4}) admits the Min-Max principle (refer to \cite{Ch,CH} also):
\begin{eqnarray}
\Upsilon_1(\Om)=\inf_{u,\Delta u\in L^2(\Om),
\int_{\Om}udx=0}\frac{\int_{\Om}|\Delta u|^2dx}{\int_{\Om}u^2dx}.
\label{Min-Max}
\end{eqnarray}
Especially, for $\Om=B_R(0)$, we have
\begin{eqnarray}
\Upsilon_1(B_R(0))=\inf_{u,\Delta u\in L^2(B_R(0)),
\int_{B_R(0)}udx=0}\frac{\int_{B_R(0)}|\Delta
u|^2dx}{\int_{B_R(0)}u^2dx}. \label{Min-Max1}
\end{eqnarray}

For $u\in  C^{\infty}(\overline{B_R(0)})$ with
$\int_{B_R(0)}u(x)dx=0$, we have the following expansion
\begin{eqnarray}
u&=&\sum_{k=1}^{\infty} \psi_k(r)g_k(\omega), \label{eigene}
\end{eqnarray}
where $(r,\om)$ stands for the polar coordinates of ${\bf R}^n$ and
\[ \psi_k(r)=\int_{S^{n-1}}
u(r,\omega)g_k(\omega)d\omega \in C^{\infty}(0,R],~~~k\ge1.\] The
following is a basic lemma, which determines the first nonzero
Neumann eigenvalues of the bi-harmonic operator $\Delta^2$ on a ball
in ${\bf R}^n$.
\begin{Lemma} \label{lemma1}
The first nonzero Neumann eigenvalue of the bi-harmonic operator
$\Delta^2$ on the ball $B_R(0)$ in ${\bf R}^n$ ($n\ge2$) is:
$\Upsilon_1(B_R(0))=\mu_1^2(R)$, where $\mu_1(R)$ is the first
nonzero Neumann eigenvalue of the Laplacian $\Delta$ on $B_R(0)$.
\end{Lemma}
\noindent{\bf Proof}. Let $J(r)=j_k(r)$ be the Bessel-type function
solving
\begin{eqnarray} r^2 J^{\prime\prime} +(n-1)r J^{\prime}
+[r^2+\lambda_k]J=0, ~~J^{\prime}(R)=0, \label{Bess}
\end{eqnarray}
and $\{\nu_{k,l}\}_{l=1}^{\infty}$ with $\nu_{k,1}<\nu_{k,2}<\cdots$
be the non-negative zeroes of $J^{'}(r)=j^{'}_k(r)$. It is
well-known that the first nonzero Neumann eigenvalue (\ref{N}) of
the Laplacian $\Delta$ on the ball $B_R(0)$ is
$\mu_1(B_{R}(0))=\nu_{1,1}$ (for example, see \cite{Wei,SY}). Due to
\cite{watson} Chapter 6 (see \cite{korenev} also), for every $k$ we
have the following Bessel-Fourier expansion,
 \[ \psi_k(r)=\sum_{l\ge 1} c_{k,l} j_k (\nu_{k,l}r),\]
where
$c_{k,l}=\int^R_0\psi(r)j_k(\nu_{k,l}r)dr/\int^R_0j^2_k(\nu_{k,l}r)dr$.
Combining this with (\ref{eigene}), we obtain that, for $u\in
C^{\infty}(\overline {B_R(0)})$,
 \[u(x)=\sum_{k\ge1,l\ge 1} c_{k,l} j_k (\nu_{k,l}r)g_k(\omega) \]
 and hence
 \[ \int_{B_R(0)}|u(x)|^2dx=\sum_{k\ge1,l\ge 1} c_{k,l}^2\int_0^R r^{n-1} |j_k(\nu_{k,l}r)|^2dr.\]
Furthermore, for $u\in  C^{\infty}(\overline {B_R(0)})$ with
$\int_{B_R(0)}u(x)dx=0$ and
 $\partial_r u=\partial_r \Delta u=0$ on $r=R$,
 we have, for any fixed $k\ge1$ and $l\ge1$,
 \begin{eqnarray}
 \int_{B_R(0)}j_k(\nu_{k,l}r)g_k(\omega)\Delta u(x) dx &=& \int_{B_R(0)}
 \Delta \bigg(j_k(\nu_{k,l}r)g_k(\omega)\bigg)u(x) dx\nonumber
  \\ &=& \int_{B_R(0)}L_k \bigg(j_k(\nu_{k,l}r)\bigg)g_k(\omega)u(x) dx\nonumber
   \\ &=&\nu_{k,l} \int_{B_R(0)} j_k(\nu_{k,l}r)g_k(\omega) u(x)
   dx, \label{iden1}
 \end{eqnarray}
 where
\[ L_k=\frac{\partial^2}{\partial r^2}+\frac{n-1}{r}\frac{\partial}{\partial r}
+\frac{\lambda_k}{r^2}.\] Here in the last equality of (\ref{iden1})
we have used the following identity
\begin{eqnarray*} \Delta \bigg(j_k(\nu_{k,l}r)g_k(\omega)\bigg)&=&\bigg(\frac{\partial^2}
{\partial r^2}+\frac{n-1}{r}\frac{\partial}{\partial r} +\frac{1}{r^2}\Delta_{S^{n-1}}
\bigg)\bigg(j_k(\nu_{k,l}r)g_k(\omega)\bigg)\\
&=&g_k(\omega) L_k \bigg(j_k(\nu_{k,l}r)\bigg) \\
&=&\nu_{k,l}j_k(\nu_{k,l}r)g_k(\omega).
\end{eqnarray*}
Hence (\ref{iden1}) leads to, for any $u\in
C^{\infty}(\overline{B_R(0)})$ with $\int_{B_1(0)}u(x)dx=0$,
\begin{eqnarray} \int_{B_R(0)} |\Delta u|^2 dx&=&\sum_{k\ge 1,l\ge
1}c_{k,l}^2(\nu_{k,l})^{2}
\int_0^R r^{n-1}|j_k(\nu_{k,l}r)|^2dr\nonumber\\
&\ge&\nu_{1,1}^2\int_{B_R(0)} |u|^2 dx\label{iden2}\end{eqnarray}
and the equality holds if and only if
$u(x)=j_1(\nu_{1,1}r)\frac{x_i}{r}, i=1,\cdots,n$. By the Min-Max
principle (\ref{Min-Max}), we see that
$$
\Upsilon_1(B_R(0))=\nu_{1,1}^2=\mu_1^2(R).$$ This proves Lemma
\ref{lemma1}. ~~$\Box$

By the way, from the proof of Lemma \ref{lemma1} we also have
\begin{Corollary} \label{Coro1}
The first (nonzero) Neumann eigenfunctions of the Laplacian $\Delta$
on the ball $B_R(0)$:
$$j_1(\nu_{1,1}r)\frac{x_i}{r},~~
i=1,\cdots,n,$$ where $j_1(r)$ solves (\ref{Bess}) with
$\lambda_k=\lambda_1=-(n-1)$, are first (nonzero) Neumann
eigenfunctions of the bi-harmonic operator $\Delta^2$ on $B_R(0)$.
\end{Corollary}

Up to the author's knowledge, besides universal inequalities on
eigenvalues of the multi-Laplacian operators (see \cite{CY,WX} and
the references there in), there are no other results in literature.
In order to understand the isoperimetric property of eigenvalues of
the multi-Laplacian operators, we may consider the following Neumann
eigenvalue problem of the even-multi-Laplacian operator
$\Delta^{2m}$ ($m\ge1$) on a bounded domain $\Om$ in ${\bf R}^n$:
\begin{eqnarray}
\left\{\begin{array}{c} \Delta^{2m}u={\hat\Upsilon} u,
~~~~~~~~~~~~~~~~~
~~~~~~~~~\hbox{in}~~~\Om~\\
\frac{\pa u}{\pa n}=\frac{\pa \Delta u}{\pa n}=\cdots=\frac{\pa
\Delta^{2m-1} u}{\pa n}=0,~\hbox{on} ~~ \pa\Om,
\end{array}\right.\label{multi}
\end{eqnarray}
where $\frac{\pa}{\pa n}$ denotes the outward normal derivative on
the boundary $\pa\Om$ of $\Om$. One notes that, when $m=1$, Problem
(\ref{multi}) is reduces to Problem (\ref{4}). It is obvious that
Problem (\ref{multi}) has discrete eigenvalues
$$0={\hat\Upsilon}_0<{\hat\Upsilon}_1\le{\hat\Upsilon}_2
\le\cdots\le {\hat\Upsilon}_k\le \cdots\nearrow \infty.$$ Now we
similarly establish the following lemma.

\begin{Lemma} \label{lemma2}
The first nonzero Neumann  eigenvalue of the even-multi-Laplacian
operator $\Delta^{2m}$ on the ball $B_R(0)$ in ${\bf R}^n$ ($n\ge2$)
is: ${\hat\Upsilon}_1(B_R(0))=\mu_1^{2m}(R)$, where $\mu_1(R)$ is
the first nonzero Neumann eigenvalue of the Laplacian $\Delta$ on
$B_R(0)$.
\end{Lemma}
\noindent{\bf Proof}. We go along the line of the proof of Lemma
\ref{lemma1}. Notice that, for $u\in C^{\infty}(\overline{B_R(0)})$
with $\int_{{B_R(0)}}u(x)dx=0$ and
 $\partial_r u=\partial \Delta u=\cdots=\partial_r \Delta^{m-1} u=0$ on $r=R$,
we have the expansion
$$u(x)=\sum_{k\ge1,l\ge 1} c_{k,l} j_k (\nu_{k,l}r)g_k(\omega)$$
and, for any fixed $k\ge1$ and $l\ge1$,
 \begin{eqnarray*}
 \int_{B_R(0)}j_k(\nu_{k,l}r)g_k(\omega)\Delta^m u(x) dx
   =(\nu_{k,l})^m \int_{B_R(0)} j_k(\nu_{k,l}r)g_k(\omega) u(x) dx.
 \end{eqnarray*}
This leads to \begin{eqnarray} \int_{B_R(0)} |\Delta^m u|^2
dx&=&\sum_{k\ge 1,l\ge 0}c_{k,l}^2(\nu_{k,l})^{2m} \int_0^R
r^{n-1}|j_k(\nu_{k,l}r)|^2dr
\nonumber\\
&\ge&(\nu_{1,1})^{2m}\int_{B_R(0)} |u|^2 dx,
\label{iden3}\end{eqnarray} and the equality holds if and only if
$u(x)=j_1(\nu_{1,1}r)\frac{x_i}{r}, i=1,\cdots,n$. Now by the
Min-Max principle for the first eigenvalue
${\hat\Upsilon}_1(B_R(0))$ of Problem (\ref{multi}):
\begin{eqnarray}
{\hat\Upsilon}_1(B_R(0))=\inf_{u,\Delta^mu \in L^{2}(B_R(0)),
\int_{B_R(0))}udx=0}\frac{\int_{B_R(0)}|\Delta
u|^{2m}dx}{\int_{B_R(0)}u^2dx},  \nonumber
\end{eqnarray}
we have
$$
{\hat\Upsilon}_1(B_R(0))=\nu_{1,1}^{2m}=\mu_1^{2m}(R).$$ This proves
Lemma \ref{lemma2}. ~~$\Box$

\begin{Corollary}
The first (nonzero) Neumann eigenfunctions of the Laplacian $\Delta$
on the ball $B_R(0)$:
$$j_1(\nu_{1,1}r)\frac{x_i}{r},~~
i=1,\cdots,n,$$ where $j_1(r)$ solves (\ref{Bess}) with
$\lambda_k=\lambda_1=-(n-1)$, are first (nonzero) Neumann
eigenfunctions of the bi-harmonic operator $\Delta^{2m}$ on
$B_R(0)$.
\end{Corollary}

Based on the basic lemma \ref{lemma1} or lemma \ref{lemma2}, the
proof of the isoperimetric inequality (\ref{?}) of
Szeg\"o-Weinberger type for the first nonzero Neumann eigenvalue of
$\Delta^2$ or $\Delta^{2m}$ becomes now clear and natural. This will
be done in the next section.

\section{\large Isoperimetric inequality for the Neumann eigenvalue of the
bi-harmonic operator} In this section, we shall first prove that the
isoperimetric inequality (\ref{?}) of Szeg\"o-Weinberger type holds
still for the first nonzero Neumann eigenvalue $\Upsilon_1(\Om)$ of
the bi-harmonic operator $\Delta^2$. Then we generalize it to the
even-multi-Laplacian operator $\Delta^{2m}$ on a bounded $\Om$ in
${\bf R}^n$.

\begin{Theorem}\label{Theorem1}
Let $\Om$ be a connected bounded smooth domain in ${\bf R}^n$
($n\ge2$), then the first nonzero Neumann eigenvalue
$\Upsilon_1(\Om)$ of the bi-harmonic operator $\Delta^{2}$ satisfies
the following inequality of Szeg\"o-Weinberger type:
\begin{eqnarray}
\Upsilon_1(\Om)\le \Upsilon_1(B_{\Om}),\label{Neuiso1}
\end{eqnarray}
where $B_{\Om}$ is a ball in ${\bf R}^n$ with the same volume of
$\Om$. Moreover, the equality holds if and only if $\Om$ is
$B_{\Om}$.
\end{Theorem}
\noindent{\bf Proof}. By the Min-Max principle, one has
$$
\Upsilon_1(\Om)=\inf_{u,\Delta u\in
L^2(\Om),\int_{\Om}u=0}\frac{\int_{\Om}|\Delta
u|^{2}dx}{\int_{\Om}|u|^2dx}.
$$
We try to construct some suitable trial functions $u$ with
$\int_{\Om}u=0$ to give an optimal upper bound for
$\Upsilon_1(\Om)$. From the basic lemma \ref{lemma1} in the previous
section, we know that
$$
g(r)\frac{x_i}{r}, ~i=1,\cdots,n
$$
are eigenfunctions corresponding to $\Upsilon_1(B_{\Om})$ of the
bi-harmonic operator $\Delta^2$  on the ball $B_{\Om}$, where $g(r)$
satisfies
\begin{eqnarray}
\left\{\begin{array}{c}\frac{d^2g}{dr^2}+\frac{n-1}{r}\frac{d
g}{dr}+\left(\mu_1-\frac{n-1}{r^2}\right)g=0\\
\frac{d g}{dr}(R)=0,~~~~~~~~~~~~~~~~~~~~~~~~~~~\end{array}\right.
\label{est0}
\end{eqnarray}
where $\mu_1$ is the first nonzero Neumann eigenvalue of $\Delta$ on
$B_{\Om}$ and $R$ is the radius of $B_{\Om}$, i.e. $B_{\Om}=B_R(0)$.
If we let $\phi(r)=r^{\frac{n-2}{2}}g(r)$, then we transform
Eq.(\ref{est0}) into
$$
\phi^{''}+\frac{1}{r}\phi^{'}+\left(\mu_1-\frac{\left(\frac{n}{2}\right)^2}
{r^2}\right)\phi=0,$$ which is a standard Bessel equation:
$$
\varphi^{''}+\frac{1}{r}\varphi^{'}+\left(1-\frac{\left(\frac{n}{2}\right)^2}
{r^2}\right)\varphi=0,$$ by a rescaling $r\to \sqrt{\mu_1}~ r$. This
implies that $g(r)$ can be explicitly given by
$$g(r)=r^{-\frac{n-2}{2}}J_{\frac{n}{2}}(\sqrt{\mu_1}~r),$$ where
$J_{\frac{n}{2}}(r)$ is the $\frac{n}{2}$-Bessel function of first
kind and $\sqrt{\mu_1}$ is the first positive zero of $g(r)$.
Now we define an auxiliary function
$$
G(r)=g(r),~~r\ge0.
$$
For any $x_0\in$ the convex hull of $\Om$, we define a vector by
$$
V(x_0)=\sum_{i=1}^n\left(\int_{\Om}\frac{(x-x_0)G(r(x,x_0))}{r(x,x_0)}\right)\frac{\pa}{\pa
x_i},$$ which can be regarded as a continuous vector field on the
convex hull of $\Om$. Furthermore, it is easy to see that $V(x_0)$
points inward at the boundary of the convex hull of $\Om$.
Therefore, it follows from the Hopf theorem that there exists an
$x^0$ in the convex hull of $\Om$ such that $V(x^0)=0$ (the detailed
arguments are referred to \cite{Wei}). By a translation in ${\bf
R}^n$, we may assume that $x^0=0$ without the loss of generality and
thus set $u_i=G(r)\frac{x_i}{r}, i=1,\cdots,n$. Then we have
$$\int_{\Om}u_i=0,~~i=1,\cdots,n.$$
We will see below that $u_i,\Delta u_i\in L^2(\Om)$ ($i=1,\cdots,n$)
by our above construction of $G(r)$. Therefore, for the bounded
smooth domain $\Om$ in ${\bf R}^n$,  by the Min-Max principle we
have
\begin{eqnarray}
\Upsilon_1(\Om)\int_{\Om}\sum_{i=1}^nu^2_i\le
\int_{\Om}\sum_{i=1}^n|\Delta u_i|^2.\label{est11}
\end{eqnarray}
One also notes that
\begin{eqnarray}
\sum_i|\Delta u_i|^2&=&\sum_{i=1}^n\left|\left(\frac{d^2}{dr^2}
+\frac{n-1}{r}\frac{d}{dr}+\frac{\Delta_{{\bf
S}^{n-1}(1)}}{r^2}\right)u_i\right|^2\nonumber\\
&=&\sum_{i=1}^n\left|\left[\left(\frac{d^2}{dr^2}
+\frac{n-1}{r}\frac{d}{dr}-\frac{n-1}{r^2}\right)G(r)\right]
\frac{x_i}{r}\right|^2\nonumber\\
&=&\left[\left(\frac{d^2}{dr^2}
+\frac{n-1}{r}\frac{d}{dr}-\frac{n-1}{r^2}\right)G(r)\right]^2,
\nonumber
\end{eqnarray}
and $\sum_iu^2_i=G^2$. Substituting them into (\ref{est11}) we have
\begin{eqnarray}
\Upsilon_1(\Om)\le \frac{\int_{\Om}\left[\left(\frac{d^2}{dr^2}
+\frac{n-1}{r}\frac{d}{dr}-\frac{n-1}{r^2}\right)G(r)\right]^2
}{\int_{\Om}G^2}. \label{est31}
\end{eqnarray}
Since $G(r)=g(r)$ and $g(r)$ solves (\ref{est0}), we see that
$$
\bigg[\left(\frac{d^2}{dr^2}
+\frac{n-1}{r}\frac{d}{dr}-\frac{n-1}{r^2}\right)G(r)\bigg]^2=\mu_1^{2}G^2(r).$$
Substituting this identity into (\ref{est31}) we finally have
\begin{eqnarray}
\Upsilon_1(\Om)\le \frac{\int_{\Om}\left[\left(\frac{d^2}{dr^2}
+\frac{n-1}{r}\frac{d}{dr}-\frac{n-1}{r^2}\right)G(r)\right]^2
}{\int_{\Om}G^2}=\frac{\mu_1^{2}\int_{\Om} G^2
}{\int_{\Om}G^2}=\mu_1^{2}=\Upsilon_1(B_{\Om}). \nonumber
\end{eqnarray}
When the equality holds, all the inequalities in the above steps
become equalities, which implies that the $n$ trial functions
$u_i(x)$ ($1\le i\le n$) should be eigenfunctions corresponding to
$\Upsilon_1(\Om)$. This leads to that $\Om$ must be the ball
$B_{\Om}$. The proof of Theorem \ref{Theorem1} is completed. $\Box$

We remark that our method used in the proof of Theorem
\ref{Theorem1} is applicable to the even-multi-Laplacian operator
$\Delta^{2m}$. That says, we similarly have

\begin{Theorem} \label{Theorem2}
Let $\Om$ be a connected bounded smooth domain in ${\bf R}^n$
($n\ge2$), then the first nonzero Neumann eigenvalue ${\hat
\Upsilon}_1(\Om)$ of the multi-Laplacian operator $\Delta^{2m}$
satisfies the following isoperimetric inequality of
Szeg\"o-Weinberger type:
\begin{eqnarray}
{\hat \Upsilon}_1(\Om)\le {\hat \Upsilon}_1(B_{\Om}).\label{Neuiso}
\end{eqnarray}
Moreover, the equality holds if and only if $\Om$ is $B_{\Om}$.
\end{Theorem}
\noindent{\bf Proof}. For the bounded smooth domain $\Om$ in ${\bf
R}^n$, by the Min-Max principle, one has
\begin{eqnarray}
{\hat\Upsilon}_1(\Om)=\inf_{f,\Delta^m f\in
L^{2}(\Om),\int_{\Om}f=0}\frac{\int_{\Om}|\Delta^m
f|^{2}}{\int_{\Om}|f|^2}. \label{Min-Maxm}
\end{eqnarray}
Now we try to construct some suitable trial functions $f$ to give an
optimal upper bound for ${\hat\Upsilon}_1(\Om)$. From the basic
lemma 2 in the previous section, we know that
$$
g(r)\frac{x_i}{r}, ~i=1,\cdots,n
$$
are first eigenfunctions of $\Delta^{2m}$ on $B_{\Om}$, where
$g(r)=r^{-\frac{n-2}{2}}J_{\frac{n}{2}}(\sqrt{\mu_1}~r)$ satisfies
(\ref{est0}), in which $J_{\frac{n}{2}}(r)$ is the
$\frac{n}{2}$-Bessel function of first kind. With the same argument
done in the proof of Theorem \ref{Theorem1}, up to  a translation in
${\bf R}^n$, we may have
$$\int_{\Om}f_i=0,~~i=1,\cdots,n,$$
where $f_i=g(r)\frac{x_i}{r}$, $i=1,\cdots,n$. It is also obvious
that $f_i,\Delta^m f_i\in L^{2}(\Om)$ ($i=1,\cdots,n$).  By the
Min-Max principle (\ref{Min-Maxm}) we see that
\begin{eqnarray}
\mu_1(\Om)\int_{\Om}\sum_{i=1}^nf^2_i\le
\int_{\Om}\sum_{i=1}^n|\Delta^m f_i|^2.\label{est1}
\end{eqnarray}
Now
\begin{eqnarray}
\sum_i|\Delta^m f_i|^2&=&\sum_{i=1}^n\left|\left(\frac{d^2}{dr^2}
+\frac{n-1}{r}\frac{d}{dr}+\frac{\Delta_{{\bf
S}^{n-1}(1)}}{r^2}\right)^mf_i\right|^2\nonumber\\
&=&\sum_{i=1}^n\left|\left[\left(\frac{d^2}{dr^2}
+\frac{n-1}{r}\frac{d}{dr}-\frac{n-1}{r^2}\right)^mg(r)\right]
\frac{x_i}{r}\right|^2\nonumber\\
&=&\left[\left(\frac{d^2}{dr^2}
+\frac{n-1}{r}\frac{d}{dr}-\frac{n-1}{r^2}\right)^mg(r)\right]^2,
\nonumber
\end{eqnarray}
and $\sum_if^2_i=g^2$. Substituting these identities into
(\ref{est1})
 we have
\begin{eqnarray}
{\hat\Upsilon}_1(\Om)\le
\frac{\int_{\Om}\left[\left(\frac{d^2}{dr^2}
+\frac{n-1}{r}\frac{d}{dr}-\frac{n-1}{r^2}\right)^mg(r)\right]^2
}{\int_{\Om}g^2}. \label{est3}
\end{eqnarray}
Since $g(r)$ satisfies (\ref{est0}), we have that the integrand
\begin{eqnarray} Q(r)&=&\left[\left(\frac{d^2}{dr^2}
+\frac{n-1}{r}\frac{d}{dr}-\frac{n-1}{r^2}\right)^mg(r)\right]^2\nonumber\\
&=&(-\mu_1)^{2m}g^2(r)=\mu_1^{2m}g^2(r). \label{es}
\end{eqnarray}
Substituting (\ref{es}) into (\ref{est3}) we finally have
\begin{eqnarray}
{\hat\Upsilon}_1(\Om)\le
\frac{\int_{\Om}\left[\left(\frac{d^2}{dr^2}
+\frac{n-1}{r}\frac{d}{dr}-\frac{n-1}{r^2}\right)g(r)\right]^2
}{\int_{\Om}g^2}=\frac{\mu_1^{2m}\int_{\Om} g^2
}{\int_{\Om}g^2}=\mu_1^{2m}={\hat \Upsilon}_1(B_{\Om}). \nonumber
\end{eqnarray}
With the same argument in the proof of Theorem \ref{Theorem1}, we
also obtain $\Om=B_{\Om}$ when the equality holds. This completes
the proof of Theorem \ref{Theorem2}. $\Box$

\bigskip
Theorem \ref{Theorem1} indicates that not only the eigenvalue
problem (\ref{4}) is actually the Neumann eigenvalue problem for the
bi-harmonic operator $\Delta^2$ on a bounded smooth doamain $\Om$ in
${\bf R}^n$, but also the isoperimetric inequality of
Szeg\"o-Weinberger type for the first nonzero eigenvalue
$\Upsilon_1(\Om)$ holds still for the bi-harmonic operator. Theorem
\ref{Theorem2} shows that the similar result is true for the
even-multi-Laplacian $\Delta^{2m}$ ($m\ge1$). This might reveal a
universal isoperimetric property of the first nonzero Neumann
eigenvalue of the multi-Laplacian $\Delta^{m}$ ($m\ge2$).

Finally we point out that whether the isoperimetric inequality
(\ref{?}) holds for the first nonzero Neumann eigenvalue of the
odd-multi-Laplacian $\Delta^{2m+1}$ is unknown at present time,
since our method applied here does not work directly in this case.
However, we conjecture that it is true. Another interesting problem
is: Whether all the first eigenfunctions of the bi-harmonic operator
$\Delta^2$ on the ball $B_{\Om}$ are just these listed in Corollary
\ref{Coro1} or not. These problems deserve for future study.

\section*{ Acknowledgement}  The authors are
supported by the National Natural Science Foundation of China (Q.D.
by Grant Nos.10971030, F.G. by J0730103 and Y.Z. by Grant No.
11031001).


\begin{thebibliography}{50}
\setlength{\itemsep}{-3pt} \small

\bibitem{A}
M.S. Ashbaugh,Isoperimetric and universal inequalities for
eigenvalues, in {\it Spectral Theory and Geometry}, London
Mathematical Society Lecture Note Series {\bf 273}, ed. by B. Davies
and Y. Safarev, Cambridge University Press, 1999, pp.95-139.

\bibitem{AB}
M.S. Ashbaugh and R.D. Benguria, Proof of the Payne-Polya-Weiberger
conjecture, Bull. Amer. Math. Soc. {\bf 25} (1991) 19-29.

\bibitem{AB1}
M.S. Ashbaugh and R.D. Benguria, A sharp bound for the ratio of the
first two eigenvalues of Dirichlet Laplacian and extensions, Ann.
Math. {\bf 135} (1992) 601-628.

\bibitem{AB2}
M.S. Ashbaugh and R.D. Benguria, On Rayleigh's conjecture for the
clamped palte and its generalization to three dimensions, Duke Math.
J. {\bf 78} (1995) 1-17.

\bibitem{Ch}
I. Chavel, Eigenvalues in Riemannian Geometry,Academic Press, Inc.,
Orlando, San Diego, New York, 1984.

\bibitem{Ch1}
I. Chavel, Isoperimetric ineqalities (differential geometric and
analytic perspectives), Cambridge Tracts in Mathematics 145,
Cambridge Univ. Press, 2001.

\bibitem{CY}
Q.M. Cheng and H.C. Yang, Universal bounds for eigenvalues of a
buckling problem, Comm. Math. Phys. {\bf 262} (2006) 663-675.

\bibitem{CH}
R. Courant and D. Hilbert, Methods of mathematical physics, Vol. I,
Wily (Intersicence), New York, 1953; Vol.II, 1967.

\bibitem{Fa}
G. Faber, Beweis, dass unter allen homogenen Membranen von gleicher
Fl\"ache und gleicher Spannung die kriesf\"ormige den tiefsten
Grundton gibt, Sitzungberichte der mathematisch-physikalischen
Klasse der Bayerschen Akademie der Wissenschaften zu M\"unchen
Jahrgabgm 1923, pp.169-172.

\bibitem{Kr}
E. Krahn, \"Uber eine von Rayleigh formulierte minimaleigenschaft
des kreises, Math. Ann. {\bf 94} (1925) 97-100.

\bibitem{korenev} B.G. Korenev, { Bessel Functions and their Applications,}
Chapman \& Hall/CRC 2002.


\bibitem{Na}
N.S. Nadirahvili, Rayleigh's conjecture on the principle frequency
of the clamped plate, Arch. Rational Mech. Anal. {\bf 129} (1995)
1-10.

\bibitem{PPW}
L.E. Payne, G. Polya and H.F. Weinberger, On the ratio of
consecutive eigenvalues, J. Math. and Phys. {\bf 35} (1956) 289-298.

\bibitem{PS} G. Polya and G. Szeg\"o, Isoperimetric inequality in
mathematical physics, Annals of Mathematics Studies, {\bf 27}
Princeton University Press, Princeton, New York, 1951.

\bibitem{Ra}
 J.W.S. Rayleigh, The theory of sound, Dover Publications, New York
1945 (republication of the 1894/96 edition).

\bibitem{SY}
R. Scheon and S.T. Yau, Lectures on differential geometry,
International Press, 1994.

\bibitem{Se}
G. Szeg\"o, Inequalities for certain eigenvalue of a membrane of
given area, J. Rational Mech. Anal. {\bf 3} (1954) 343-356.

\bibitem{WX}
Q.L. Wang and C.Y. Xia, Universal bounds for eigenvalues of the
biharmonic operator on Riemannian manifolds, J. Funct. Anal. {\bf
245} (2007) 334-352.

\bibitem{watson} G. N. A. Watson,  A Treatise on the theory of Bessel
functions, 2 ed. Cambridge, the University Press 1959.


\bibitem{Wei}
H.F. Weinberger, An isoperimetric inequality for the $n$-dimensional
free membrane problem, J. Rational Mech. Anal. {\bf 5} (1956)
633-636.

\end{thebibliography}
\end{document}